\documentclass[11pt,oneside]{article}
\usepackage[T1]{fontenc}
\usepackage[utf8]{inputenc}
\usepackage{amsmath,amssymb,amsthm,mathtools}
\usepackage{newpxtext,newpxmath}
\usepackage[paperwidth=190.5mm,paperheight=254mm,left=19mm,right=19mm,top=19mm,bottom=20mm,headsep=6mm]{geometry}
\usepackage{microtype,needspace,booktabs,array,enumitem}
\usepackage{xcolor,xurl,fancyhdr}
\usepackage[colorlinks=true,linkcolor=blue!35!black,citecolor=blue!35!black,urlcolor=blue!35!black]{hyperref}
\usepackage{bookmark}
\hypersetup{pdftitle={An Explicit Construction of Commuting Quantum KdV Hamiltonians in the Universal Virasoro Vertex Algebra},pdfsubject={Symmetric Virasoro PBW action, reduced polynomial inverse, quantum conservation and commutativity}}
\newcommand{\Q}{\mathbb Q}
\newcommand{\C}{\mathbb C}
\newcommand{\base}{\mathbf k}
\newcommand{\Vcal}{\mathcal V}
\newcommand{\Pcal}{\mathcal P}
\newcommand{\Fcal}{\mathcal F}
\newcommand{\wt}{\operatorname{wt}}
\newcommand{\ad}{\operatorname{ad}}
\newcommand{\Res}{\operatorname{Res}}
\newcommand{\Span}{\operatorname{span}}
\newcommand{\vac}{\mathbf 1}
\newcommand{\XH}{X_{\mathrm H}}
\newcommand{\hh}{\hbar}
\newcommand{\eps}{\epsilon}
\newcommand{\ev}{\pi}
\newcommand{\rep}{\rho}
\newcommand{\extract}{\mathcal C}
\newcommand{\mon}{\mathfrak m}
\newtheorem{theorem}{Theorem}[section]
\newtheorem{proposition}[theorem]{Proposition}
\newtheorem{lemma}[theorem]{Lemma}

\theoremstyle{definition}
\newtheorem{definition}[theorem]{Definition}

\theoremstyle{remark}

\makeatletter
\renewenvironment{proof}[1][\proofname]{\par
  \pushQED{\qed}\normalfont\topsep6\p@\@plus6\p@\relax
  \trivlist\item[\hskip\labelsep\bfseries\itshape #1\@addpunct{.}]\ignorespaces
}{\popQED\endtrivlist\@endpefalse}
\makeatother
\numberwithin{equation}{section}
\allowdisplaybreaks[2]
\setlist{nosep,leftmargin=1.6em}

\title{\bfseries An Explicit Construction of\\Commuting Quantum KdV Hamiltonians in the Universal Virasoro Vertex Algebra}
\author{Hongbo Zhao}
\date{4 November 2026}
\begin{document}
\maketitle
\vspace{-6mm}

\begin{abstract}
We give an explicit formula for elements $H_{2k}$ of conformal weight
$2k$, $k\ge1$, in the universal Virasoro vertex operator algebra $V$
at every central charge, such that $H_{2k}(0)H_{2l}\in TV$ for all
$k,l\ge1$. Modulo translation, the kernel of $H_4(0)$ is generated
by $[H_{2k}]$ in each positive even weight and is zero in positive
odd weights. The formula uses a realization on differential polynomials
by possibly infinite-order differential operators in a completion of
the Weyl algebra, together with the inverse of the identity plus a
nilpotent operator on each weight space. A Koszul-type exactness
argument for the classical KdV flows proves conservation; uniqueness
and the Jacobi identity give commutativity. The classical limit recovers
the normalized Gelfand--Dickey densities. Explicit low-weight elements
and exact symbolic verifications illustrate the construction.

\end{abstract}
\begingroup\footnotesize\tableofcontents\endgroup
\newpage

\section{Introduction}\label{sec:intro}
The Korteweg--de Vries equation $u_t=u_{xxx}+6uu_x$ admits an
infinite family of conserved quantities that commute under the KdV
Poisson bracket. Their Hamiltonian flows form the KdV hierarchy
\cite{Lax,GD}. The first densities are
\[
 u,\qquad u^2,\qquad u^3-\tfrac12u_x^2.
\]
Their integrals are constant for smooth periodic or rapidly decaying
solutions. The corresponding algebraic statement is a local
conservation law: the time derivative of a density is a total spatial
derivative. Thus the equation and its conservation laws can be formulated entirely
in the algebra of differential polynomials
$\Pcal_0:=\Q[u,u_{-1},u_{-2},\ldots]$, where
$u_{-j}:=\partial_x^ju$ and the total derivative satisfies
$\partial u_{-j}=u_{-j-1}$. The KdV evolution acts on differential
polynomials by a derivation $X_0$, determined by
$X_0u=u_{-3}+6uu_{-1}$ and $X_0\partial=\partial X_0$.
Explicitly,
\[
 X_0:=\sum_{j\ge0}\partial^j(u_{-3}+6uu_{-1})\partial_{u_{-j}}.
\]
The quotient $\Pcal_0/\partial\Pcal_0$ is the space of local
functionals; it expresses integration by parts without choosing
boundary conditions. For the second KdV Poisson bracket, $X_0$ is
the Hamiltonian flow generated by $[u^2]$
(Section~\ref{sec:classical-hamiltonians}).
A class $[f]$ is a conserved quantity for the
flow $X_0$ precisely when
\[
 X_0[f]=0,\qquad\text{equivalently}\qquad X_0f\in\partial\Pcal_0.
\]
Pairwise commutativity is a further property of the chosen conserved
quantities. The pseudodifferential formulation places KdV within the
KP hierarchy \cite{DSKV}.

The Virasoro vertex operator algebra provides the quantum version of
this problem \cite{FF,BLZ}. Let $V^c$ be the universal Virasoro VOA of central charge
$c$, with conformal vector $\omega$ and translation $T$. We use modes
$Y(a,z)=\sum_{n\in\mathbb Z}a(n)z^{-n-1}$; in particular,
\begin{equation}\label{eq:omegaField}
 Y(\omega,z)=\sum_{n\in\mathbb Z}L_nz^{-n-2},\qquad
 \omega(n)=L_{n-1}.
\end{equation}
The space $V^c/TV^c$
is a Lie algebra with bracket
\[
 \{[a],[b]\}:=[a(0)b].
\]
Our result concerns this universal VOA at every central charge,
including $c=0$.
\begin{quote}
\textbf{Main result (Theorem~\ref{thm:main}).}
There is an explicit formula for vectors $H_{2k}^{(c)}\in V^c$,
of weight $2k$ and with highest PBW term $\omega(-1)^k\vac$, such that
\[
 H_{2k}^{(c)}(0)H_{2l}^{(c)}\in TV^c\qquad(k,l\ge1).
\]
The kernel of $(\omega(-1)\omega)(0)$ on $V^c/TV^c$ is spanned
by $[H_{2k}^{(c)}]$ in each positive even weight and is zero in
positive odd weights.
\end{quote}
Thus the conserved classes for the weight-four action are unique
up to scalar in each positive even weight; the displayed identity
also asserts their pairwise commutativity.
The formula is \eqref{eq:Fk}, followed by the change of variables
\eqref{eq:reverse}; the first vectors are displayed in
Appendix~\ref{app:computations}.

Earlier explicit low-weight densities were obtained by Sasaki and
Yamanaka \cite{SY}. Feigin and Frenkel establish existence and
commutativity of quantum integrals of motion for generic parameter
by a homological deformation argument
\cite[Section~4.7]{FF}. Bazhanov, Lukyanov, and Zamolodchikov
construct quantum monodromy and commuting transfer operators, and
relate their asymptotic expansions to local integrals \cite{BLZ}.
They also display the first densities, of weights two, four, and six
\cite[equations~(5)--(6)]{BLZ}. More recent explicit expressions,
including weights eight and ten, appear in
\cite[Appendix~D]{MNRT}. These approaches provide existence
and spectral information; here the aim is a direct formula for the
polynomial densities in fixed VOA coordinates. The proof works directly with identities in the universal Virasoro
VOA, independently of transfer-operator asymptotics. Our local-functional
bracket uses the zeroth products $a(0)b$ in the mode convention above.
Comparison with cylinder currents
also requires the exponential change of coordinates. In Zhu's
square-bracket formalism \cite{Zhu}, this is expressed by
\[
 Y[a,z]:=Y(e^{zL_0}a,e^z-1),\qquad
 \widetilde\omega:=\omega-\frac{c}{24}\vac.
\]
Thus comparison with the formulas of \cite{BLZ,MNRT} requires
converting both coordinates and normal-ordering conventions.
We do not need this change of coordinates in the construction below.

For the construction, we introduce a parameter $\hh$ and pass to
polynomial coordinates obtained by symmetrization. The quantum KdV action then becomes
an operator $X_\hh$ deforming $X_0$, and the conserved densities take
the form
\begin{equation}\label{eq:introformula}
 \boxed{F_k=(1+M)^{-1}u^k.}
\end{equation}
Here $M$ acts nilpotently at each weight after reduction modulo total
derivatives, so the inverse is a finite sum. Sections~\ref{sec:vir}--\ref{sec:formula}
give the precise setup; equation~\eqref{eq:reverse} returns to the
original central-charge coordinates.

Four steps explain the construction and its proof.
\begin{enumerate}[label=(\roman*)]
\item Symmetric PBW coordinates realize the Virasoro modes as possibly
infinite-order differential operators, finite on each polynomial.
Their coefficients have degree at most one in the coordinates.
\item The weight-four action has a part that can be inverted explicitly
after reduction modulo total derivatives. The remaining correction
strictly lowers a finite ordering on each weight space. This reduces
the construction to inverting the identity plus a nilpotent operator,
giving \eqref{eq:introformula}.
\item The commuting classical third and fifth KdV flows define a
Koszul-type complex. Its middle exactness, together with the explicit weight-six
conservation identity, shows that the constructed $F_k$ are
conserved densities for $X_\hh$.
\item The same inverse determines the kernel in each weight. Odd-weight
vanishing and the Jacobi identity then imply pairwise commutativity.
\end{enumerate}

The paper is organized as follows. Section~\ref{sec:vir} develops the
differential representation, and
Section~\ref{sec:weightfour} computes the weight-four action.
Section~\ref{sec:reduction} constructs reduction and extraction;
Section~\ref{sec:formula} gives the formula and states the main theorem.
Section~\ref{sec:classical}
proves classical exactness, and Section~\ref{sec:quantum} completes
the conservation, uniqueness, and specialization arguments. The
appendix records explicit elements and computational checks. Readers
primarily interested in the formula may start with
\eqref{eq:Fk}, Theorem~\ref{thm:main}, and Appendix~\ref{app:computations}.
The arguments for the differential realization, finite inverse,
and quantum conservation and commutativity are given explicitly.
We use standard PBW and VOA identities, with conventions as in
\cite{Kac,LL}, and recall the classical Gelfand--Dickey construction
where it enters the proof \cite{DSKV}.

\section{Symmetric PBW coordinates and differential operators}\label{sec:vir}

\subsection{Virasoro modes and symmetric PBW coordinates}
Introduce a parameter $\hh$ and work over the
polynomial ring
\[
 \base:=\Q[\hh].
\]
The completion $\Q[[\hh]]$ is used only in the lifting argument of
Section~\ref{sec:quantum}. After adjoining $\hh^{-1}$, set
$U:=\hh\omega$ and $c:=6/\hh$; this rescales the Virasoro relations
so that all commutators vanish at $\hh=0$ while their first-order
terms retain the classical bracket. Let $\Vcal$ be the universal
vertex algebra over $\base$ generated by $U$ with singular products
\begin{equation}\label{eq:regularIntro}
 U(z)U(w)\sim\hh\left(
 \frac{3}{(z-w)^4}+\frac{2U(w)}{(z-w)^2}
                       +\frac{\partial_wU(w)}{z-w}\right).
\end{equation}
Here $U(z):=Y(U,z)$, and $\sim$ denotes equality of singular parts.
Equivalently, its modes satisfy the algebraic relation
\begin{equation}\label{eq:UmodeIntro}
 [U(m),U(n)]=\hh(m-n)U(m+n-1)
       +\frac\hh2m(m-1)(m-2)\delta_{m+n,2}.
\end{equation}
The family $\Vcal$ is a graded vertex algebra over $\base$ and has
conformal vector $U/\hh$ after adjoining $\hh^{-1}$. Its
$\hh=0$ specialization is the classical differential algebra.

The mode convention \eqref{eq:omegaField}
and the scaling $U=\hh\omega$ motivate the following indices:
nonpositive indices will correspond directly to the differential
polynomial generators $u_{-j}$. Define
\begin{equation}\label{eq:Kmodes}
 U_{-j}:=j!U(-j-1)\quad(j\ge0),\qquad
 U_n:=U(n-1)\quad(n\ge1),\qquad U_n\vac=0\quad(n\ge1).
\end{equation}
Thus $U_{-j}\vac=T^jU$, matching $u_{-j}=\partial^ju$ without
additional factorials in the symmetrization map. Positive indices
label annihilation modes. After localization,
\begin{equation}\label{eq:scaling}
 \begin{gathered}
 U=\hh\omega,\qquad c=6/\hh,\qquad
 U_{-j}=j!\hh L_{-j-2}\ (j\ge0),\quad U_n=\hh L_{n-2}\ (n\ge1),\\
 [L_m,L_n]=(m-n)L_{m+n}
              +\frac c{12}(m^3-m)\delta_{m+n,0}.
 \end{gathered}
\end{equation}
By \eqref{eq:Kmodes}, $U_n=\max(0,-n)!\,U(n-1)$; hence
\eqref{eq:UmodeIntro}, with indices $m-1,n-1$, gives for all $m,n\in\mathbb Z$
\[
 \begin{aligned}
 [U_m,U_n]
 &=\hh(m-n)\frac{\max(0,-m)!\,\max(0,-n)!}
                    {\max(0,2-m-n)!}U_{m+n-2}\\
 &\quad+\frac\hh2\max(0,-m)!\,\max(0,-n)!
        \bigl((m-2)^3-(m-2)\bigr)\delta_{m+n,4}.
 \end{aligned}
\]
Here $\max(0,-n)!$ is $(-n)!$ for $n\le0$ and $1$ for $n\ge1$.
The induced vacuum module carries the vertex algebra $\Vcal$.
Its translation satisfies $T\vac=0$ and
\[
 [T,U_{-j}]=U_{-j-1}\quad(j\ge0),\qquad
 [T,U_n]=-(n-1)U_{n-1}\quad(n\ge1).
\]
PBW gives the basis $U_{-r_1}\cdots U_{-r_m}\vac$ with
$r_1\ge\cdots\ge r_m\ge0$. Standard field reconstruction gives
the vertex-algebra structure \cite{Kac}. At fixed $c$, $V^c$ is the universal Virasoro VOA. Its vacuum module
is the highest-weight-zero Verma module modulo the submodule generated
by $L_{-1}\vac$. We work with this universal algebra rather than
its simple quotient.

The vectors $T^jU$ suggest the differential polynomial coordinates
\begin{equation}\label{eq:jets}
 \Pcal:=\base[u_{0},u_{-1},u_{-2},\ldots],\qquad
 u:=u_{0},\qquad\partial:=\sum_{j\ge0}u_{-j-1}\partial_{u_{-j}}.
\end{equation}
Thus $\Pcal=\base\otimes_{\Q}\Pcal_0$, and specializing $\hh=0$
recovers the algebra $\Pcal_0$ introduced above. We use capitals for
the generating vectors and modes on the Virasoro side, and lowercase
letters for their commuting polynomial coordinates.
Jet order always means derivative order: $u_{-j}$ has order $j$, with subscript $-j$. Two gradings will be used:
\begin{enumerate}[label=(\roman*)]
\item the conformal weight, $\wt(u_{-j}):=j+2$;
\item the \emph{factor number}, $\deg_{\mathrm f}(u_{-j}):=1$, extended
additively to products.
\end{enumerate}
Thus $uu_{-2}$ and $u_{-1}^2$ have factor number two and weight six,
whereas $u^3$ has factor number three and the same weight. The
parameter $\hh$ has both degrees zero. Let $\Pcal^{[m]}$ denote the
factor-homogeneous subspace of degree $m$, and put
$\Pcal^{\le m}:=\bigoplus_{r=0}^m\Pcal^{[r]}$.

Since $U_{-j}\vac=T^jU$ for $j\ge0$, define the \emph{symmetrization map}
\begin{equation}\label{eq:sym}
 \sigma(u_{-j_1}\cdots u_{-j_m})
   :=\frac1{m!}\sum_{\tau\in S_m}
       U_{-j_{\tau(1)}}\cdots U_{-j_{\tau(m)}}\vac,
 \qquad\sigma(1):=\vac.
\end{equation}
The average includes repeated indices. This is the standard PBW
symmetrization underlying the symbol calculus of Berezin
\cite{Berezin}; the induced-module version needed here is proved below.

\begin{proposition}\label{prop:sym}
The map $\sigma:\Pcal\to\Vcal$ is a weight-preserving filtered
$\base$-linear isomorphism, and
\begin{equation}\label{eq:symT}
 \sigma\partial=T\sigma.
\end{equation}
\end{proposition}
\begin{proof}
Let $\tau:\Pcal\to\Vcal$ send each monomial to its ordered word in
the $U_{-j}$. PBW makes $\tau$ an isomorphism. Filter by factor number.
Exchanging adjacent creation modes introduces a commutator and hence
a word with one fewer factor. Consequently,
\[
 (\tau^{-1}\sigma-1)\Pcal^{\le m}\subseteq\Pcal^{\le m-1}.
\]
For $Q:=\tau^{-1}\sigma-1$, this gives $Q^{m+1}=0$ on
$\Pcal^{\le m}$ and the inverse
$\sigma^{-1}=(\sum_{r=0}^m(-Q)^r)\tau^{-1}$. This proves the
filtered isomorphism even with infinitely many generators.
Finally $[T,U_{-j}]=U_{-j-1}$; translation acts on each factor of a
symmetrized word in turn, giving \eqref{eq:symT}.
\end{proof}
The map $\sigma$ is linear; it does not identify polynomial
multiplication with $a(-1)b$. For example,
\begin{equation}\label{eq:symexample}
 \sigma(uu_{-1})=(U_{-1}U_0+\tfrac\hh{12}U_{-3})\vac,
 \qquad\sigma(u_{-3})=U_{-3}\vac.
\end{equation}

\subsection{Symbols of operators on differential polynomials}
The first step is to express the Virasoro action as operators on
differential polynomials. Transport the mode action by symmetrization:
\begin{equation}\label{eq:Ln}
 \rep(U_n):=\sigma^{-1}U_n\sigma\in\operatorname{End}_{\base}(\Pcal).
\end{equation}
Every such endomorphism can be written uniquely as a possibly infinite
sum of differential operators with polynomial coefficients. We now
specify the formal series that encode this expression; no bound on
the differential order is assumed.

The auxiliary commuting parameters $t_j$ and multi-index entries are indexed by derivative order $j\ge0$, so $t_j$ pairs with $u_{-j}$ and $U_{-j}$. A multi-index $I=(I_i)$ has nonnegative integer entries, all but
finitely many of which are zero. We use
\[
 u^I:=\prod_i u_{-i}^{I_i},\qquad I!:=\prod_i I_i!,\qquad
 |I|:=\sum_i I_i,\qquad \partial_u^I:=\prod_i\partial_{u_{-i}}^{I_i}.
\]
For a $\base$-module $W$, put
\[
 W[[t]]:=\prod_I Wt^I,\qquad [t^I]\Bigl(\sum_Jw_Jt^J\Bigr):=w_I.
\]
Thus $[t^I]$ extracts the coefficient of the monomial $t^I$.
We give this space the product topology in which $W$ is discrete.
If $W$ is an algebra, each product coefficient is a finite sum over
$J\le I$. Power series in an element of positive total $t$-degree are
therefore well defined. A $\base$-linear map extends coefficientwise
to a continuous $\base[[t]]$-linear map. The formal derivatives
$\partial_{t_i}$ are $\base$-linear, but not $\base[[t]]$-linear;
they commute with these coefficientwise extensions.

Set $U(t):=\sum_{j\ge0}t_jU_{-j}$. Here $U_{-j}$ is a creation
operator, whereas $U_{-j}\vac=T^jU$ is a vector. In the centrally extended
mode algebra $\mathfrak g$ specified below,
\[
 U(t)\in\mathfrak g[[t]],\qquad
 e^{U(t)}\in U(\mathfrak g)[[t]],\qquad
 e^{U(t)}\vac\in\Vcal[[t]].
\]
On the polynomial side, let $u(t):=\sum_j u_{-j}t_j$ and
$e^{u(t)}:=\sum_I u^It^I/I!\in\Pcal[[t]]$.
The coefficientwise extension of $\sigma$ satisfies
\begin{equation}\label{eq:symexponential}
 \boxed{\sigma(e^{u(t)})=e^{U(t)}\vac.}
\end{equation}
Indeed, each coefficient on either side is the same average of
creation words. Thus symmetric PBW corresponds to the exponential of
their sum. No completion in the $u$ variables is needed: every
$t$-coefficient is a polynomial. The right side is module-valued;
we invert exponentials only in the respective algebras.

For a $\base$-linear endomorphism $D$ of $\Pcal$, define
\begin{equation}\label{eq:normal}
 \operatorname{symb}(D)(u,t)
   :=e^{-u(t)}D(e^{u(t)})
    =\sum_I c_I(u)t^I\in\Pcal[[t]].
\end{equation}
Comparing coefficients after multiplication by $e^{u(t)}$ gives
\begin{equation}\label{eq:normaloperator}
 \frac{D(u^J)}{J!}
  =\sum_{I\le J}c_I(u)\frac{u^{J-I}}{(J-I)!},\qquad
 D=\sum_I c_I(u)\partial_u^I.
\end{equation}
Equivalently,
\[
 c_I(u)=\sum_{J\le I}
       \frac{(-u)^{I-J}}{(I-J)!\,J!}D(u^J).
\]
All coefficient multiplications stand to the left of the derivatives.
Only finitely many derivatives contribute on a polynomial, although
their orders need not be uniformly bounded. This reconstruction holds
for every linear operator; the coordinate-degree bound for a mode
will follow from its explicit symbol. The sums in
\eqref{eq:normaloperator} form the completion of the polynomial Weyl
algebra used here: they act on each polynomial by a finite sum.

\subsection{Explicit mode symbols from the exponential map}
Berezin's symbol calculus gives differential formulas for the regular
enveloping-algebra action \cite{Berezin}. The induced-module
exponential formula used below is a specialization of
\cite[Theorem~23 and equations~(4.4)--(4.5)]{Molotkov}.
We include a coefficientwise proof adapted to the present Virasoro
module and its countably many modes.
The preceding reconstruction applies to any endomorphism. For a Lie
algebra action we need a more explicit symbol, in particular a bound
on the coordinate degree of its coefficients. This bound will isolate
the terms that preserve or raise factor number in Section~\ref{sec:weightfour}.
To retain the central term before evaluating it on the vacuum, let
$\mathfrak g$ have the mode relations above with their scalar central
term multiplied by $\mathbf c/6$. The central generator $\mathbf c$
acts as the scalar $6$ on the vacuum, so the displayed scalar mode
relations are recovered. It is distinct from the central charge $c$.
We use the splitting into Lie subalgebras
\[
 \mathfrak g=\mathfrak g_-\oplus\mathfrak g_+,\qquad
 \mathfrak g_-:=\Span_{\base}\{U_{-j}:j\ge0\},\quad
 \mathfrak g_+:=\Span_{\base}\{U_n:n>0\}\oplus\base\mathbf c,
\]
with projections denoted by $-$ and $+$. More explicitly,
\[
 \Vcal=U(\mathfrak g)\otimes_{U(\mathfrak g_+)}\base\vac,
\]
where $U_n$ for $n>0$ annihilates $\vac$ and $\mathbf c$ acts as $6$.
PBW multiplication gives
$U(\mathfrak g)\simeq U(\mathfrak g_-)\otimes_{\base}U(\mathfrak g_+)$;
after imposing the vacuum action, symmetrization identifies
$\Vcal\simeq U(\mathfrak g_-)\simeq\Pcal$ as $\base$-modules.

Set
\begin{equation}\label{eq:dexpfunctions}
 \varphi(s):=\frac{1-e^{-s}}s,\qquad
 \psi(s):=\varphi(s)^{-1}=\frac{s}{1-e^{-s}}
       =1+\frac{s}{2}+\frac{s^2}{12}-\frac{s^4}{720}+\cdots.
\end{equation}
These quotients denote power series. For $X\in\mathfrak g$, define
\begin{equation}\label{eq:mode-directions}
 \boxed{\begin{aligned}
 v_X(U,t)&:=\psi(\ad U(t))\bigl(e^{-\ad U(t)}X\bigr)_-
                    =\sum_{j\ge0}v_{X,j}(t)U_{-j},\\
 \bigl(e^{-\ad U(t)}X\bigr)_+\vac&=6c_X(t)\vac,
 \end{aligned}}
\end{equation}
where $c_X(t)$ is the coefficient of $\mathbf c$ in the conjugated
series. The projection precedes $\psi$. All scalar coefficients lie
in $\base[[t]]$, so their dependence on $\hh$ is implicit.
For each $t^I$, only its finitely many input indices and at most
$|I|$ brackets contribute. Hence every coefficient of $v_X(U,t)$
is a finite linear combination of the $U_{-j}$, and the shorthand
\[
 v_X(u,t)\vac:=\sigma^{-1}(v_X(U,t)\vac)
                  =\sum_{j\ge0}v_{X,j}(t)u_{-j}
\]
is well defined in $\Pcal[[t]]$, since $\sigma(u_{-j})=U_{-j}\vac$.
The next identity is the formal version of the differential of the
exponential map; see \cite[Section~5.4, Theorem~5.4]{Hall}.
\begin{lemma}[Formal differentiation of the exponential]\label{lem:dexp}
In $U(\mathfrak g)[[t]]$, the following identity holds coefficientwise:
\begin{equation}\label{eq:dexpseries}
 \partial_{t_j}e^{U(t)}
 =\int_0^1 e^{(1-s)U(t)}U_{-j} e^{sU(t)}\,ds
 =e^{U(t)}\varphi(\ad U(t))U_{-j}.
\end{equation}
\end{lemma}
\begin{proof}
To verify this identity formally, expand the exponential and differentiate
each power: for $p,q\ge0$, the word $U(t)^pU_{-j}U(t)^q$ has coefficient
$1/(p+q+1)!$. Expanding the integrand gives the same coefficient,
since the integral of $(1-s)^ps^q/(p!q!)$ over $[0,1]$ is
$1/(p+q+1)!$. Integration here is performed coefficientwise on
polynomials in $s$, so no analytic convergence is required.
Factoring out $e^{U(t)}$ and using conjugation gives the last equality.
\end{proof}

For arbitrary $X\in\mathfrak g$, extend \eqref{eq:Ln} by
$\rep(X):=\sigma^{-1}X\sigma$.
\begin{proposition}[Symbol of a Lie algebra action]\label{prop:affine}
For $X\in\mathfrak g$, the symbol and action of $X$ are
\begin{equation}\label{eq:affine}
 \boxed{\begin{aligned}
 \ell_X(u,t)&:=\operatorname{symb}(\rep(X))(u,t)
                 =v_X(u,t)\vac+6c_X(t),\\
 \rep(X)&=\sum_I\bigl([t^I]\ell_X(u,t)\bigr)\partial_u^I.
 \end{aligned}}
\end{equation}
In particular, every coefficient has coordinate degree at most one.
\end{proposition}

\begin{proof}
\textbf{Step 1. Move the action to the vacuum.}
Conjugation in the completed enveloping algebra gives
\[
 Xe^{U(t)}=e^{U(t)}\bigl(e^{-\ad U(t)}X\bigr).
\]
Each iterated commutator lies in $\mathfrak g[[t]]$ and raises total
$t$-degree. Its coefficient at any fixed $t^I$ is a finite linear
combination of modes and the central generator. The same therefore
holds for formal power series in $\ad U(t)$ applied to $X$.
Splitting into the two subalgebras and acting on the vacuum gives
\begin{equation}\label{eq:modeconjugation}
 Xe^{U(t)}\vac
 =e^{U(t)}\bigl(e^{-\ad U(t)}X\bigr)_-\vac
       +6c_X(t)e^{U(t)}\vac.
\end{equation}
Since $\mathfrak g_-$ is a subalgebra, $\ad U(t)$ preserves it and
raises total $t$-degree. Thus $\varphi(\ad U(t))$ is invertible there,
and the definition of $v_X$ says exactly that
\begin{equation}\label{eq:directionequation}
 \varphi(\ad U(t))v_X(U,t)
       =\bigl(e^{-\ad U(t)}X\bigr)_-.
\end{equation}

\textbf{Step 2. Transfer by symmetrization.}
Apply Lemma~\ref{lem:dexp} to express derivatives of the noncommutative
exponential, then transfer them to polynomial coordinates.
The extension of $\sigma$ acts on polynomial coefficients only, so
$\partial_{t_j}\sigma(\sum_I f_It^I)
 =\sum_I\sigma(f_I)\partial_{t_j}t^I
 =\sigma(\partial_{t_j}\sum_I f_It^I)$.
Differentiating the PBW exponential therefore gives
\begin{equation}\label{eq:PBWdirection}
 \partial_{t_j}(e^{U(t)}\vac)
       =\sigma(\partial_{t_j}e^{u(t)})
       =\sigma(u_{-j} e^{u(t)}).
\end{equation}
Multiply these two derivative identities by $v_{X,j}(t)$ and sum:
\begin{equation}\label{eq:modePBWdirection}
 \begin{aligned}
 e^{U(t)}\varphi(\ad U(t))v_X(U,t)\vac
 &=\sum_{j\ge0}v_{X,j}(t)\partial_{t_j}(e^{U(t)}\vac)\\
 &=\sigma\!\left[v_X(u,t)\vac\,e^{u(t)}\right].
 \end{aligned}
\end{equation}
The sum is coefficientwise finite. The scalar coefficients multiply
after differentiation and pass through the $\base[[t]]$-linear map
$\sigma$; no derivative acts on them.

\textbf{Step 3. Recover the symbol.} Substitution into the conjugation
identity now yields
\begin{equation}\label{eq:modepreimage}
 Xe^{U(t)}\vac
  =\sigma\!\left[(v_X(u,t)\vac+6c_X(t))e^{u(t)}\right].
\end{equation}
Apply $\sigma^{-1}$ and cancel $e^{u(t)}$ in the symbol
definition. This proves \eqref{eq:affine}, including the degree bound.
\end{proof}

For the Virasoro modes needed below, take $X=U_n$ and write
$\ell_n:=\ell_{U_n}$ and $c_n:=c_{U_n}$. Thus \eqref{eq:affine} gives
\[
 \ell_n(u,t)=v_{U_n}(u,t)\vac+6c_n(t),\qquad
 \rep(U_n)=\sum_I\bigl([t^I]\ell_n(u,t)\bigr)\partial_u^I.
\]
For creation modes $n\le0$, all brackets stay in $\mathfrak g_-$,
so $c_n=0$. To illustrate an annihilation mode, set $t_0=t$ and
$t_j=0$ for $j\ge1$. Then $U(t)=tU_0$ and
\[
 e^{-t\ad U_0}U_4
  =U_4+4\hh tU_2+\frac{\hh t}{2}\mathbf c
                         +4\hh^2t^2U_0.
\]
The creation component commutes with $U(t)$, so
$v_{U_4}(U,t)=4\hh^2t^2U_0$ and $6c_4(t)=3\hh t$. Consequently,
\begin{equation}\label{eq:K2restriction}
 \ell_4(u;t,0,\ldots)=3\hh t+4\hh^2ut^2,\qquad
 \rep(U_4)|_{\base[u]}=3\hh\partial_u+4\hh^2u\partial_u^2.
\end{equation}
The second-order term is invisible on $1$ and $u$; coordinate degree
at most one does not mean differential order at most one.

\Needspace{7\baselineskip}
\begin{lemma}[Parameter degree of a mode]\label{lem:modeparameter}
For a multi-index $I$, every coefficient of $[t^I]\ell_n$ is
homogeneous of degree $|I|$ in $\hh$.
\end{lemma}
\begin{proof}
Each commutator with $U(t)$ contributes one parameter and one factor
of $\hh$. A term of degree $|I|$ in \eqref{eq:mode-directions} has
$|I|$ brackets in total, before and after projection; the central
coefficient $c_n$ obeys the same count. Projection, factorials in
$U_{-j}=j!U(-j-1)$, and evaluation of $\mathbf c$ introduce only rational
factors. Equation~\eqref{eq:affine} proves the assertion.
\end{proof}

\section{Leading mode terms and the weight-four zero mode}\label{sec:weightfour}
The construction depends on the zero-mode action of
$H_4=U(-1)U$. We first isolate the terms of the individual modes
with the largest factor-number shifts, then combine them to determine
this action. Here ``leading'' refers to factor number, not to
differential order or degree in $\hh$.
\subsection{Mode terms that raise or preserve factor number}
For a differential term $p(u)\partial_u^I$ with factor-homogeneous coefficient $p$, the factor-number shift is $\deg_{\mathrm f}p-|I|$.
Write $D^{[s]}$ for the component of an operator $D$ that changes
factor number by $s$, so $D^{[s]}\Pcal^{[m]}\subseteq\Pcal^{[m+s]}$.
By \eqref{eq:affine}, every single-mode coefficient has degree at
most one. Hence terms of differential order at least two lower
factor number. The multiplication and first-order coefficients of
any $D$ are recovered by
\[
 c_{\boldsymbol 0}=D(1),\qquad c_{\boldsymbol e_j}=D(u_{-j})-u_{-j}D(1).
\]
Here $\boldsymbol 0$ is the zero multi-index and $\boldsymbol e_j$
has its only nonzero entry, equal to one, at position $j$.
This determines the following pieces exactly, without determining
the higher-order terms.

For $i\ge0$, $U_{-i}\vac=\sigma(u_{-i})$.
The multiplication coefficient of $\rep(U_{-i}):=\sigma^{-1}U_{-i}\sigma$
is $u_{-i}$. Symmetric ordering gives
\[
 U_{-i}U_{-j}\vac-\sigma(u_{-i} u_{-j})=\tfrac12[U_{-i},U_{-j}]\vac,
 \qquad [U_{-i},U_{-j}]=\hh\frac{(j-i)i!j!}{(i+j+2)!}U_{-i-j-2}.
\]
Consequently,
\begin{equation}\label{eq:creation}
 \begin{aligned}
 \rep(U_{-i})&=u_{-i}
 +\frac\hh2\sum_{j\ge0}\frac{(j-i)i!j!}{(i+j+2)!}
                     u_{-i-j-2}\partial_{u_{-j}}\\
 &\quad+\text{terms of differential order at least two}.
 \end{aligned}
\end{equation}
The two displayed pieces have factor shifts $+1$ and $0$; all omitted
terms lower factor number.

For $n\ge-1$, $\rep(U_{n+2})(1)=0$, and the first-order coefficient is
$\sigma^{-1}([U_{n+2},U_{-j}]\vac)$. On the vacuum, the commutator is
\[
 [U_{n+2},U_{-j}]\vac=
 \begin{cases}
 \displaystyle\hh\frac{(n+j+2)j!}{(j-n)!}U_{-j+n}\vac,&j\ge n,\\
 0,&j<n
 \end{cases}
 \; +\frac\hh2(n^3-n)j!\delta_{n,j+2}\vac.
\]
Thus the factor-preserving part is
\begin{equation}\label{eq:annihilation}
 \boxed{\rep(U_{n+2})^{[0]}=
 \hh\sum_{j\ge\max(0,n)}
       \frac{(n+j+2)j!}{(j-n)!}\,u_{-j+n}\partial_{u_{-j}},
 \qquad n\ge-1.}
\end{equation}
For $n\ge2$, the remaining first-order term is the scalar operator
\begin{equation}\label{eq:centralfirst}
 \boxed{\frac\hh2(n^3-n)(n-2)!\,\partial_{u_{-n+2}}.}
\end{equation}
It lowers factor number by one. Higher-order terms also lower factor
number, so \eqref{eq:centralfirst} is not the complete shift-$(-1)$
component. In particular, the second-order term in
\eqref{eq:K2restriction} must not be discarded from the full action.
For the first two modes the following formulas are exact:
\begin{equation}\label{eq:translationgrading}
 \rep(U_1)=\hh\partial,\qquad
 \rep(U_2)=\hh\sum_{j\ge0}(j+2)u_{-j}\partial_{u_{-j}}.
\end{equation}
Indeed $U_1=\hh T$ on the vacuum module, and $U_2$ acts by
$\hh$ times the weight. Their commutators with each creation factor
preserve symmetric ordering, so no higher-order correction occurs.

\subsection{The weight-four action}
Define
\begin{equation}\label{eq:H4def}
 H_4:=U(-1)U,\qquad
 \boxed{X_\hh:=\sigma^{-1}\hh^{-1}H_4(0)\sigma.}
\end{equation}
Because the $(-1)$-product is normal ordered, its two summands at
mode zero coincide after reindexing. Thus
\begin{equation}\label{eq:Xmodes}
 \boxed{X_\hh=\frac2\hh\sum_{r\ge2}\frac1{(r-2)!}\rep(U_{-r+2})\rep(U_{r-1}).}
\end{equation}
On an input of weight $w$, terms with $r>w+3$ vanish because the
right mode lowers weight by $r-3$. Every such right mode has index
at least $1$ and is divisible by $\hh$: it kills the vacuum, and
every commutation through a creation word supplies a factor of
$\hh$. The displayed division therefore defines an operator over
$\base$, not an operator with a pole at $\hh=0$.

For operators whose right symbol has coordinate degree at most one,
the Leibniz rule gives
\[
 \operatorname{symb}(D_1D_2)
 =\operatorname{symb}(D_1)\operatorname{symb}(D_2)
 +\sum_j(\partial_{t_j}\operatorname{symb}(D_1))
          (\partial_{u_{-j}}\operatorname{symb}(D_2)).
\]
Derivatives of order at least two kill the right coefficients.
At a fixed coefficient in $t$, those coefficients involve finitely
many coordinates. Applied to \eqref{eq:Xmodes}, this yields
\begin{equation}\label{eq:Xsymbol}
 \operatorname{symb}(X_\hh)=\frac2\hh\sum_{r\ge2}\frac1{(r-2)!}\left(
       \ell_{2-r}\ell_{r-1}
       +\sum_{j\ge0}(\partial_{t_j}\ell_{2-r})
                         (\partial_{u_{-j}}\ell_{r-1})\right).
\end{equation}
In particular, its coordinate degree is at most two.

\begin{proposition}\label{prop:Xstructure}
The operator $X_\hh$ commutes with $\partial$, raises conformal
weight by three, and raises factor number by at most one. Its
factor-increasing part is $\XH$, where
\begin{equation}\label{eq:Hopf}
 \XH:=6\sum_{j\ge0}\partial^j(uu_{-1})\partial_{u_{-j}}.
\end{equation}
Consequently,
\begin{equation}\label{eq:R}
 \boxed{R:=X_\hh-\XH}
\end{equation}
preserves or lowers factor number, and all its coefficients lie in
$\Q[\hh]$.
\end{proposition}
\begin{proof}
A zero-mode operator commutes with translation. Together with
$\sigma\partial=T\sigma$, this gives $[X_\hh,\partial]=0$.
Each product of modes in \eqref{eq:Xmodes} raises weight by three.
The creation factor raises factor number by at most one, and the
annihilation factor does not raise it.

More precisely, the shift-$+1$ part can only arise from the
multiplication term in \eqref{eq:creation} composed with the
factor-preserving part \eqref{eq:annihilation}. All other combinations
have shift at most zero. The resulting operator is first-order, and
its coefficient at $\partial_{u_{-j}}$ is
\begin{equation}\label{eq:Hopfcoefficient}
 \begin{aligned}
 (X_\hh)^{[+1]}u_{-j}
 &=2j!\sum_{a=0}^{j+1}
       \frac{j+a+1}{a!(j+1-a)!}\,u_{-a} u_{-j-1+a}\\
 &=3\partial^{j+1}(u^2)=6\partial^j(uu_{-1}).
 \end{aligned}
\end{equation}
For the middle equality, pair the terms indexed by $a$ and $j+1-a$;
the average of their numerator factors is $3(j+1)/2$. This identifies
the whole shift-$+1$ operator as \eqref{eq:Hopf}, not merely its value
on $u$. Polynomial dependence on $\hh$ was established before
\eqref{eq:Xsymbol}.
\end{proof}

\subsection{The classical limit and the factor-preserving terms}
By Lemma~\ref{lem:modeparameter}, a differential-order-$b$ coefficient
of a single mode carries $\hh^b$. In \eqref{eq:Xmodes}, after division
by $\hh$, the only terms surviving at $\hh=0$ are creation
multiplication times the \emph{complete} first-order annihilation
part. The differentiation terms in the composition formula require a
positive-order contribution from each mode, hence vanish after division
by $\hh$ and specialization. The coordinate-linear part gives $\XH$.
With $n=r-3\ge2$,
the scalar part \eqref{eq:centralfirst} contributes
\[
 \frac2\hh\frac{u_{-n-1}}{(n+1)!}
 \frac\hh2(n^3-n)(n-2)!\,\partial_{u_{-n+2}}
       =u_{-n-1}\partial_{u_{-n+2}}.
\]
Hence
\begin{equation}\label{eq:Xclassicalmodes}
 \boxed{X_\hh\big|_{\hh=0}
   =\XH+\sum_{j\ge0}u_{-j-3}\partial_{u_{-j}}=X_0.}
\end{equation}
This explains directly why the dispersive term survives: it comes
from a central term that lowers the factor number of a single mode.
The PBW-leading part and the classical limit are different operations.

\textbf{The action on the first two powers of $u$.}
On the one-generator vector only $r=2,3,5$ in \eqref{eq:Xmodes}
contribute. Before applying $\sigma^{-1}$ the result is
\[
 (2U_{0}U_{-1}+4U_{-1}U_{0}+U_{-3})\vac.
\]
Using \eqref{eq:symexample} gives
\begin{equation}\label{eq:Xu}
 \boxed{X_\hh u=6uu_{-1}+\left(1-\frac\hh6\right)u_{-3}.}
\end{equation}
Applying the same mode sum to $u^2$ and reordering its terms gives
the following identity. The ordered-mode calculation and the
symmetrization conversions are recorded in Appendix~\ref{app:modecheck}:
\begin{equation}\label{eq:Xu2}
 X_\hh u^2
 =12u^2u_{-1}+\left(2+\frac{5\hh}{3}\right)uu_{-3}
  +3\hh u_{-1}u_{-2}+\left(\frac{3\hh}{20}-\frac{\hh^2}{5}\right)u_{-5}.
\end{equation}

For a factor-homogeneous input, let $R^{[s]}$ be the part changing
factor number by $s$: $R^{[s]}\Pcal^{[m]}\subseteq\Pcal^{[m+s]}$.
Thus $R=\sum_{s\le0}R^{[s]}$. The coordinate-degree-two bound now
makes the factor-preserving action on a polynomial in $u$ alone
completely explicit:
\begin{equation}\label{eq:preserving}
 \boxed{R^{[0]}f(u)=\left(1-\frac\hh6\right)u_{-3}f'(u)
          +\hh\left(uu_{-3}+\frac32u_{-1}u_{-2}\right)f''(u).}
\end{equation}
Indeed a factor-preserving term with coefficients on the left must have equal coordinate
degree and differential order; the latter is therefore at most two.
Equation~\eqref{eq:Xu} determines the first-order coefficient. Half
the two-factor part of $X_\hh(u^2)-2uX_\hh(u)$ is
$\hh(uu_{-3}+\tfrac32u_{-1}u_{-2})$, determining the second-order coefficient.
The one-factor $u_{-5}$ term of \eqref{eq:Xu2} cannot contribute here.

These truncations have distinct purposes. Creation multiplication and
factor-preserving annihilation prove the Hopf-leading statement; the
scalar first-order annihilation terms are also needed for the classical
limit; \eqref{eq:preserving} supplies the seed calculation. None of
these truncated formulas is the full remainder $R$. Evaluating the
higher-weight construction uses the complete symbols
\eqref{eq:mode-directions} and \eqref{eq:affine}, or equivalently finite
PBW commutation, in \eqref{eq:Xmodes}.

\section{Local functionals, reduction, and extraction}\label{sec:reduction}
\subsection{Local functionals and reduced representatives}\label{sec:normalform}
A \emph{local functional} is a differential polynomial considered
modulo total derivatives: it is the algebraic counterpart of the
integral of a density. These classes form the $\base$-module
\begin{equation}\label{eq:Fspace}
 \Fcal:=\Pcal/\partial\Pcal.
\end{equation}
The image $\partial\Pcal$ need not be an ideal, so this quotient is
not an algebra for polynomial multiplication. By \eqref{eq:symT},
$\Fcal$ is identified with $\Vcal/T\Vcal$. Its Lie bracket is
\begin{equation}\label{eq:localbracket}
 \{[a],[b]\}:=[a(0)b].
\end{equation}
Indeed, $(Ta)(0)=0$ and $[a(0),T]=0$ give well-definedness;
skew-symmetry and the mode commutator identity give antisymmetry and
Jacobi modulo $T\Vcal$. The bracket has weight
$\wt(a)+\wt(b)-1$.

The choice below follows the elimination of monomials linear in their
highest derivative used for scalar polynomial densities in
\cite[Section~1.3 of the arXiv version, Step~1]{HG}.
We specify the reduction map and prove uniqueness of the resulting
representative in our polynomial setting.
\begin{definition}\label{def:reduced}
A monomial is \emph{reduced} if it is a power of $u$, including $1$,
or its highest-order jet occurs at least twice.
\end{definition}
For example, $u^4,uu_{-1}^2,u_{-2}^2$ are reduced, whereas $u^2u_{-2}$ and
$u_{-1}u_{-3}$ are not. If the highest-order jet $u_{-j}$, $j\ge2$, occurs once,
write the monomial as $u_{-j} u_{-j+1}^a f$ with
$f\in\base[u,\ldots,u_{-j+2}]$. Formal integration by parts gives
\begin{equation}\label{eq:redrule}
 u_{-j} u_{-j+1}^a f\equiv
 -\frac{u_{-j+1}^{a+1}}{a+1}\partial f\pmod{\partial\Pcal}.
\end{equation}
Every resulting term has lower maximum derivative order. Define reduction by
\begin{equation}\label{eq:projectoralgorithm}
 \Pi(\mon):=
 \begin{cases}
 \mon,&\mon\text{ reduced},\\
 0,&\mon=u^a u_{-1},\\
 -\dfrac1{a+1}\Pi(u_{-j+1}^{a+1}\partial f),
 &\mon=u_{-j} u_{-j+1}^a f,\ j\ge2,
             \ f\in\base[u,\ldots,u_{-j+2}].
 \end{cases}
\end{equation}
The last case assumes that $u_{-j}$ occurs once. Recursion terminates
because the highest derivative order decreases; the second case follows from
$u^a u_{-1}=\partial(u^{a+1})/(a+1)$. Extend $\Pi$ linearly.

\begin{proposition}[Reduced representatives]\label{prop:red}
The reduced monomials form a basis of $\Fcal$. Reduction preserves
weight and factor number, and
\begin{equation}\label{eq:redproperties}
 \Pcal=\operatorname{im}\Pi\oplus\partial\Pcal,
 \qquad\Pi^2=\Pi,\quad\Pi\partial=0,\quad[f]=[\Pi f].
\end{equation}
\end{proposition}
\begin{proof}
The terminating reduction proves spanning and preserves both gradings.
For independence, let $g$ be nonconstant and $u_{-r}$ its
highest-order jet. The part of $\partial g$ containing the next jet is
$u_{-r-1}\partial g/\partial u_{-r}$. It is nonzero and contains
$u_{-r-1}$ once; all remaining terms have lower maximum derivative order. Thus
$\partial g$ is not a linear combination of reduced monomials.
For $r=0$ the same argument uses $u_{-1}g'(u)$.
\end{proof}
We henceforth represent each class by its unique reduced polynomial.
An operator $D$ commuting with $\partial$ acts on classes by
$[f]\mapsto[Df]$; its action in the reduced model is $\Pi D$.
We use the same symbol for these induced actions when the domain is
$\Fcal$. Thus $Df=0$ in $\Fcal$ means $\Pi Df=0$, not necessarily
$Df=0$ in $\Pcal$. Maps defined only on reduced representatives,
such as the extraction map below, are not applied to arbitrary
representatives of a class.
Derivative-free evaluation is the separate map
\begin{equation}\label{eq:ev}
 \ev(f):=f(u,0,0,\ldots).
\end{equation}
It descends to $\Fcal$ because $\ev\partial=0$, and $\ev\Pi=\ev$.
Thus $\Pi$ selects a representative, while $\ev$ extracts its component
in $\base[u]$. Useful reductions are
\begin{equation}\label{eq:redexamples}
 \begin{aligned}
 \Pi(u^2u_{-2})&=-2uu_{-1}^2,&\Pi(uu_{-4})&=u_{-2}^2,\\
 \Pi(u^2u_{-3})&=u_{-1}^3,&\Pi(uu_{-1}u_{-2})&=-\tfrac12u_{-1}^3.
 \end{aligned}
\end{equation}
At odd weight, classes with one or two factors vanish. One-factor
positive derivatives are exact; for two factors,
$[u_{-i} u_{-j}]=(-1)^i[uu_{-i-j}]$, and reversing the factors makes this
class zero when $i+j$ is odd.

\subsection{A partial inverse of the dispersionless action}\label{sec:direct}
The derivation $\XH$ corresponds to the dispersionless equation
$u_t=6uu_x$, often called the Hopf equation. Its action modulo total
derivatives contains a term that is diagonal after removing one
factor $u_{-1}$. To describe this term, introduce the diagonal
degree operator and its eigenvalue
\begin{equation}\label{eq:epsilon}
 E:=u_{-1}\partial_{u_{-1}}+\sum_{j\ge2}(j+1)u_{-j}\partial_{u_{-j}},
 \qquad E\mon=\eps(\mon)\mon,
 \qquad\eps(\mon):=a_1+\sum_{j\ge2}(j+1)a_j
\end{equation}
for $\mon=\prod_j u_{-j}^{a_j}$. This is a different grading from factor
number: $E$ assigns weights $0,1,3,4,\ldots$ to $u,u_{-1},u_{-2},u_{-3},\ldots$,
whereas factor number assigns weight one to every jet.
Define the splitting derivation
\begin{equation}\label{eq:N}
 N:=6\sum_{j\ge3}\sum_{a=2}^{j-1}
           \binom ja u_{-a} u_{-j-1+a}\partial_{u_{-j}}.
\end{equation}
Expansion of $\partial^j(uu_{-1})$ on the generators gives
\begin{equation}\label{eq:Hopfdecomposition}
 \XH=6u\partial+6u_{-1}E+N.
\end{equation}
For example, the terms outside $u\partial$ in the unscaled expansion
are $u_{-1}^2$ at $u_{-1}$, and
$(j+1)u_{-1}u_{-j}+\sum_{a=2}^{j-1}\binom ja u_{-a} u_{-j-1+a}$ at $u_{-j}$
for $j\ge2$. Since $u\partial\mon\equiv-u_{-1}\mon$ modulo
$\partial\Pcal$, \eqref{eq:Hopfdecomposition} implies
\begin{equation}\label{eq:diagonal}
 \Pi\XH\mon=\Pi\bigl(6(\eps(\mon)-1)u_{-1}\mon+N\mon\bigr).
\end{equation}

The following map inverts the contribution
$6(\eps(\mon)-1)u_{-1}\mon$ in \eqref{eq:diagonal} on reduced
monomials containing derivatives.
\begin{definition}[Partial inverse of the diagonal term]\label{def:B}
For a reduced monomial $\mon$, set
\begin{equation}\label{eq:B}
 \extract(\mon):=
 \begin{cases}
 \dfrac{\mon'}{6(\eps(\mon')-1)},
 &\mon=u_{-1}\mon',\quad \mon'\text{ reduced and containing a positive-order jet},\\
 0,&\text{otherwise}.
 \end{cases}
\end{equation}
Extend $\extract$ $\base$-linearly on reduced representatives.
\end{definition}
In the nonzero case $\eps(\mon\prime)>1$, so this is division by a positive
integer. The map lowers weight by three and factor number by one;
its image contains a positive-order jet and at least two factors. For example,
\begin{equation}\label{eq:Bexamples}
 \extract(u_{-1}^3)=\tfrac16u_{-1}^2,\qquad
 \extract(u_{-1}u_{-2}^2)=\tfrac1{30}u_{-2}^2,\qquad
 \extract(u^a u_{-1}^2)=\extract(u_{-2}^2)=\extract(u^a)=0.
\end{equation}

\begin{lemma}[Extraction identity]\label{lem:direct}
On reduced representatives,
\begin{equation}\label{eq:directidentity}
 \extract\Pi\XH=1-\ev+\extract\Pi N.
\end{equation}
\end{lemma}
\begin{proof}
By linearity it suffices to consider a reduced monomial $\mon$.
Applying $\extract$ to \eqref{eq:diagonal} gives
\[
 \extract\Pi\XH\mon
 =6(\eps(\mon)-1)\extract\Pi(u_{-1}\mon)
       +\extract\Pi N\mon.
\]
If $\mon$ contains a positive-order jet, then $u_{-1}\mon$ is reduced:
its highest-order jet still occurs at least twice. Hence
$\Pi(u_{-1}\mon)=u_{-1}\mon$, and Definition~\ref{def:B} gives
\[
 6(\eps(\mon)-1)\extract\Pi(u_{-1}\mon)=\mon.
\]
Here $\eps(\mon)>1$, and $\ev(\mon)=0$, so this is exactly
$(1-\ev)\mon$.

Otherwise $\mon=u^a$ for some $a\ge0$. In this case
$u_{-1}\mon=\partial(u^{a+1})/(a+1)$, so $\Pi(u_{-1}\mon)=0$,
while $(1-\ev)\mon=0$ as well. Thus in both cases the first term
in the displayed identity equals $(1-\ev)\mon$, proving
\eqref{eq:directidentity}.
\end{proof}
For example,
\[
 N(uu_{-3}^2)=36uu_{-3}u_{-2}^2,\qquad
 \Pi N(uu_{-3}^2)=-12u_{-1}u_{-2}^3,\qquad
 \extract\Pi N(uu_{-3}^2)=-\tfrac14u_{-2}^3.
\]
The map $\extract\Pi N$ preserves factor number but lowers the
higher-jet exponents in a finite order. We prove this decrease,
together with the contribution of $R$, in the next proposition.

\section{The explicit construction and the main theorem}\label{sec:formula}
\subsection{The inverse formula on each weight space}
Write $\Fcal_w$ for the component of conformal weight $w$, represented
by reduced monomials of that weight. It is a free $\base$-module of
finite rank. On this reduced model of $\Fcal$, define
\begin{equation}\label{eq:M}
 \boxed{M:=\extract\Pi(N+R).}
\end{equation}
The decomposition $X_\hh=\XH+R$ and \eqref{eq:directidentity} give
\begin{equation}\label{eq:projectedoperator}
 \boxed{\extract\Pi X_\hh=1+M-\ev.}
\end{equation}

The parameter degree below is ordinary polynomial degree in $\hh$.
It measures the coefficients, whereas conformal weight and factor
number measure the jet monomials; $\hh$ has degree zero for the latter
two gradings.
\begin{proposition}[Finite inverse on each weight space]\label{prop:finite}
On every weight space $\Fcal_w$, the operator $M$ is nilpotent.
For $k\ge1$, the reduced polynomial
\begin{equation}\label{eq:Fk}
 \boxed{F_k:=(1+M)^{-1}u^k
     =\sum_{r=0}^{d_{2k}-1}(-M)^ru^k,\qquad
 d_w:=\operatorname{rank}_{\base}\Fcal_w}
\end{equation}
has $\ev F_k=u^k$. Every remaining term has between two and $k-1$
factors. A coefficient of a term with $m$ factors has degree at most
$k-m$ in $\hh$.
\end{proposition}
\begin{proof}
\textbf{Step 1. A finite decreasing order.}
For a reduced monomial $\mon=\prod_{j=0}^{w-2}u_{-j}^{a_j}$ of weight
$w$, put
\[
 \mu_w(\mon):=(\deg_{\mathrm f}\mon,a_{w-2},a_{w-3},\ldots,a_2),
\]
with the higher-jet list empty for $w<4$, and compare the tuples
lexicographically. At fixed weight and factor number, the exponents
of $u_{-j}$ for $j\ge2$ determine those of $u,u_{-1}$, so this orders the
finite reduced basis without ties.

A term of $N$ replaces one $u_{-j}$ by two jets of strictly smaller
derivative orders. It increases factor number by one and decreases the
higher-jet tuple. Each integration-by-parts step in
\eqref{eq:redrule} decreases that tuple again, whereas extraction
removes only $u_{-1}$. Thus $\extract\Pi N$ preserves factor number
and strictly lowers the tuple. The part $\extract\Pi R$ decreases
factor number: $R$ does not increase it, reduction preserves it,
and extraction removes one factor. Every term of $M\mon$
therefore has a smaller tuple.

List the reduced basis as $\mon_1,\ldots,\mon_{d_w}$ in increasing
order. Then $M\mon_i\in\Span\{\mon_1,\ldots,\mon_{i-1}\}$, so
$M^{d_w}=0$, proving \eqref{eq:Fk}. Also $\ev M=0$, and the image
of extraction has a positive-order jet and at least two factors. Starting
from $u^k$, the first nonzero step decreases factor number because
$Nu^k=0$. This proves the asserted range of factors.

\textbf{Step 2. The parameter degree.}
Lemma~\ref{lem:modeparameter} shows that a single-mode symbol
coefficient of differential order $b$ has
$\hh$-degree $b$. In \eqref{eq:Xsymbol}, a product term of order $b$
has parameter degree $b-1$ and coordinate degree at most two; a
contraction term has parameter degree $b$ and coordinate degree at
most one. Thus a term of coordinate degree $a$ and differential
order $b$ has parameter degree at most $1-a+b$.

For a term of $R$ the factor shift is $a-b\le0$. Extraction lowers
it once more, so a term of $\extract\Pi R$ lowering factor number
by $r\ge1$ has parameter degree at most $r$. The component
$\extract\Pi N$ preserves factor number and is independent of
$\hh$. Along a chain from $k$ factors to $m$ factors, these degrees
sum to at most $k-m$.
\end{proof}
\begin{proposition}\label{prop:projectedconservation}
For every $k\ge1$, the reduced polynomial $F_k$ defined by
\eqref{eq:Fk} satisfies
\begin{equation}\label{eq:projectedeq}
 \extract\Pi X_\hh F_k=0.
\end{equation}
\end{proposition}
\begin{proof}
By Proposition~\ref{prop:finite}, $M^{d_{2k}}=0$ on $\Fcal_{2k}$.
Multiplying the finite sum in \eqref{eq:Fk} by $1+M$ therefore gives
\[
 (1+M)F_k
 =\bigl(1-(-M)^{d_{2k}}\bigr)u^k=u^k.
\]
Every term in the image of $M$ contains a positive-order jet, since
$M=\extract\Pi(N+R)$. Thus $\ev M=0$, and applying $\ev$ to
\eqref{eq:Fk} gives $\ev F_k=u^k$. Finally, the operator identity
\eqref{eq:projectedoperator}, applied to the reduced polynomial $F_k$,
yields
\[
 \extract\Pi X_\hh F_k=(1+M)F_k-\ev F_k=u^k-u^k=0.
\]
\end{proof}
Proposition~\ref{prop:projectedconservation} gives only the extracted
conservation equation. The main theorem requires
$\Pi X_\hh F_k=0$, or equivalently $X_\hh[F_k]=0$ in $\Fcal$;
this does not require $X_\hh F_k=0$ as a polynomial.
Since $\extract$ has a kernel, it cannot simply be cancelled.
Section~\ref{sec:classical} proves an injective test for compatible
classical residual pairs. Section~\ref{sec:quantum} lifts that test
to the quantum flows and applies it to $(X_\hh F_k,Y_\hh F_k)$
for $k\ge2$, thereby removing extraction from the conservation
condition. The case $k=1$ is checked directly below.
\Needspace{10\baselineskip}
\begin{theorem}[Main theorem]\label{thm:main}
Let $F_k$ be given by \eqref{eq:Fk} and put $H_{2k}:=\sigma(F_k)$
for $k\ge1$. Then:
\begin{enumerate}[label=(\roman*)]
\item \emph{Conservation and commutativity.} The polynomials $F_k$
have weight $2k$ and satisfy
\begin{equation}\label{eq:mainconservation}
 X_\hh[F_k]=0\quad\text{in }\Fcal,
\end{equation}
and
\begin{equation}\label{eq:maincommutation}
 H_{2k}(0)H_{2l}\in T\Vcal\qquad(k,l\ge1).
\end{equation}
\item \emph{Uniqueness.} At every positive weight $w$,
\[
 \ker(X_\hh:\Fcal_w\to\Fcal_{w+3})=
 \begin{cases}\base\cdot[F_k],&w=2k,\\0,&w\text{ odd}.\end{cases}
\]
Thus the normalized reduced conserved density with derivative-free
part $u^k$ is $F_k$.
\item \emph{Classical limit.} With $L:=\partial^2+u$,
\begin{equation}\label{eq:mainclassical}
 F_k\big|_{\hh=0}=I_k,\qquad
 I_k:=\frac{4^k}{\binom{2k}{k}}\Pi\Res_\partial L^{k-1/2}.
\end{equation}
\item \emph{Every central charge.} The inverse scaling
\eqref{eq:reverse} defines polynomials $F_k^{(c)}$ over $\Q[c]$.
The vectors $H_{2k}^{(c)}:=\sigma_\omega(F_k^{(c)})$, where
$\sigma_\omega$ uses jet coordinates $v_{-j}$ for $T^j\omega$, have
commuting classes in $V^c/TV^c$ for every $c\in\C$. A coefficient
of an $m$-factor term of $F_k^{(c)}$ has degree at most $k-m$ in $c$. The centralizer of
$[\omega(-1)\omega]$ is spanned by $[H_{2k}^{(c)}]$ in weight $2k$
and is zero in positive odd weights, including at $c=0$.
\end{enumerate}
\end{theorem}

For $k=1$, the formula gives $F_1=u$ and $H_2=U$. Equation
\eqref{eq:Xu} shows $X_\hh u\in\partial\Pcal$. Since
$U(0)=\hh T$, its class is central in $\Vcal/T\Vcal$; hence its
conservation and commutativity are immediate. Also $I_1=u$ by the
residue formula. The residual argument below is needed only for
$k\ge2$.

Proposition~\ref{prop:finite} already proves the explicit formula,
its normalization, and the parameter-degree bound. The weight-six
calculation below gives a second commuting zero mode.
Section~\ref{sec:classical} proves the classical exactness needed to
remove $\extract$ from \eqref{eq:projectedeq}; Section~\ref{sec:quantum}
then proves conservation, uniqueness, commutativity, and the claimed
specializations, completing Theorem~\ref{thm:main}.

\subsection{The weight-six conserved quantity and its action}
The weight-six calculation supplies more than the first nontrivial example:
it produces a second operator $Y_\hh$ commuting with $X_\hh$.
This commutation relation will make the two conservation residuals
compatible, as required by the test in Section~\ref{sec:classical}.
For weight four, $M u^2=0$ and $F_2=u^2$.
To find the weight-six member, apply \eqref{eq:preserving} to $u^3$.
Both $X_\hh u^3$ and $X_\hh u_{-1}^2$ have weight nine. All their
terms with at most two factors vanish after reduction, by the
one- and two-factor observation following \eqref{eq:redexamples}.
The factor-preserving part of $Ru^3$ is given by
\eqref{eq:preserving}, while $Ru_{-1}^2$ has at most two factors. Hence
\[
 \Pi R u^3
   =\Pi\bigl((3+\tfrac{11}2\hh)u^2u_{-3}
                     +9\hh uu_{-1}u_{-2}\bigr)
   =(3+\hh)u_{-1}^3,
 \qquad \Pi R(u_{-1}^2)=0.
\]
Since $\extract(u_{-1}^3)=u_{-1}^2/6$ and $N$ annihilates $u^3,u_{-1}^2$, the inverse
formula gives
\begin{equation}\label{eq:F3}
 \boxed{F_3=u^3-\left(\frac12+\frac\hh6\right)u_{-1}^2.}
\end{equation}
Also $\Pi\XH(u^3)=0$ and $\Pi\XH(u_{-1}^2)=6u_{-1}^3$.
Consequently,
\begin{equation}\label{eq:seed}
 \Pi X_\hh F_3
   =(3+\hh)u_{-1}^3
     -6\left(\frac12+\frac\hh6\right)u_{-1}^3=0.
\end{equation}

Choose the following ordered weight-six vector, which represents the same local functional
as the symmetric vector $\sigma(F_3)$ modulo translation:
\begin{equation}\label{eq:J6}
 \widetilde H_6:=U(-1)(U(-1)U)
       +\left(\frac12+\frac\hh6\right)(T^2U)(-1)U.
\end{equation}
Indeed, the two terms are $U_{0}^3\vac$ and $U_{-2}U_{0}\vac$, respectively, and
\begin{equation}\label{eq:J6translation}
 \sigma(F_3)-\widetilde H_6
    =-\left(\frac12+\frac\hh6\right)
          T(U_{-1}U_{0}\vac).
\end{equation}
Both vectors therefore have the same zero-mode operator. By
\eqref{eq:seed}, $H_4(0)\widetilde H_6\in T\Vcal$, so
$[a(0),b(0)]=(a(0)b)(0)$ gives
$[H_4(0),\widetilde H_6(0)]=0$. The following operators consequently
commute after localization and, by freeness, over $\base$:
\begin{equation}\label{eq:Ydef}
 \boxed{Y_\hh:=\sigma^{-1}\frac2\hh\widetilde H_6(0)\sigma,
 \qquad [X_\hh,Y_\hh]=0.}
\end{equation}
The factor $2$ makes the coefficient of $u_{-5}$ in the classical limit
equal to one; see \eqref{eq:hlimit}.
The operator has coefficients in $\base$: the specialization at
$\hh=0$ is a commutative vertex algebra, whose nonnegative products
vanish. Hence $\widetilde H_6(0)\Vcal\subseteq\hh\Vcal$.
Its weight shift is five, and it commutes with translation, so it
also acts on $\Fcal$.
\section{Classical KdV flows and a compatibility criterion}\label{sec:classical}
The classical coefficient algebra is $\Pcal_0=\Pcal/\hh\Pcal$.
Put $\Fcal_0:=\Pcal_0/\partial\Pcal_0$. The reduced basis identifies
$\Fcal$ with $\base\otimes_\Q\Fcal_0$; write $\Fcal_{0,w}$ for its
classical weight-$w$ part. Define the evolutionary derivations on the generator $u$ by
\begin{equation}\label{eq:XY}
 \begin{aligned}
 X_0u&:=u_{-3}+6uu_{-1},\\
 Y_0u&:=u_{-5}+10uu_{-3}+20u_{-1}u_{-2}+30u^2u_{-1},
 \end{aligned}
\end{equation}
and extend uniquely by $[X_0,\partial]=[Y_0,\partial]=0$ and
the Leibniz rule. Thus $X_0u_{-j}=\partial^jX_0u$ and
$Y_0u_{-j}=\partial^jY_0u$. On $\Fcal_0$, we identify classes with
their reduced representatives, so their actions include reduction by $\Pi$.
We first identify these derivations as the classical limits of the
weight-four and weight-six VOA actions and relate them to the KdV
hierarchy. We then prove that, for a pair $(a,b)$ in the required odd
weights, the two conditions $Y_0a=X_0b$ and $\extract a=0$ force
$a=b=0$. This is the classical ingredient needed to remove $\extract$
from the conservation equation.

\subsection{Classical limits and KdV conserved quantities}\label{sec:classical-hamiltonians}
The specialization of $\Vcal$ at $\hh=0$ is commutative. Dividing
the singular products in \eqref{eq:regularIntro} by $\hh$ and taking
this specialization gives the Poisson vertex bracket
\begin{equation}\label{eq:pvabracket}
 \{u_\lambda u\}=(\partial+2\lambda)u+\tfrac12\lambda^3.
\end{equation}
Thus the constant terms of $X_\hh,Y_\hh$ are evolutionary
derivations of the commutative differential algebra. Let $J$ be the
second KdV Hamiltonian operator,
\[
 J:=\tfrac12\partial^3+2u\partial+u_{-1}.
\]
The classical symbols of $H_4,\widetilde H_6$ are $u^2$ and
$u^3+\frac12uu_{-2}$ respectively. For this classical identification only, write
the variational derivative as
\[
 \delta f:=\sum_{j\ge0}(-\partial)^j\partial_{u_{-j}}f.
\]
It vanishes on total derivatives and therefore depends only on the
class of $f$. The Hamiltonian formula
$u\mapsto J\delta f$ \cite[Section~1.4]{DSKV} then gives
\begin{equation}\label{eq:hlimit}
 \boxed{X_\hh=X_0+O(\hh),\qquad Y_\hh=Y_0+O(\hh).}
\end{equation}
Indeed
\[
 J\delta(u^2)=u_{-3}+6uu_{-1},\quad
 2J\delta(u^3+\tfrac12uu_{-2})
   =u_{-5}+10uu_{-3}+20u_{-1}u_{-2}+30u^2u_{-1}.
\]
The second expression confirms the factor $2$ in \eqref{eq:Ydef}.
The limiting first-order property follows from the commutative
specialization of the $(-1)$-product; it does not assert that the
higher $\hh$-coefficients obey polynomial Leibniz rules.

We use the classical Gelfand--Dickey construction \cite{GD}.
A formulation in the present differential-algebra language, including
both the KdV and KP hierarchies, is given in \cite[Section~3]{DSKV}. Only the third and fifth
KdV flows and their common conserved densities are needed here.

\begin{lemma}[Classical conservation and normalization]\label{lem:rho}
For $L:=\partial^2+u$, the reduced polynomials $I_k$ of
\eqref{eq:mainclassical} have weight $2k$ and satisfy
\begin{equation}\label{eq:rho}
 \ev I_k=u^k,\qquad X_0I_k=Y_0I_k=0
                 \quad\text{in }\Fcal_0.
\end{equation}
Moreover $[X_0,Y_0]=0$. In particular,
$I_1=u$, $I_2=u^2$, and $I_3=u^3-\tfrac12u_{-1}^2$.
\end{lemma}
\begin{proof}
The classical result supplies the commuting Lax flows and their
residue integrals. We record the normalization check to identify them
with \eqref{eq:XY}. Formal multiplication uses
\begin{equation}\label{eq:pdomult}
 \partial^a f=\sum_{r\ge0}\binom ar(\partial^r f)\partial^{a-r},
                     \qquad a\in\mathbb Z.
\end{equation}
The monic root $Q:=L^{1/2}=\partial+\sum_{r\ge1}q_r\partial^{-r}$
is uniquely determined by $Q^2=L$: each new coefficient occurs with
coefficient two. Operators have order bounded above, so every
coefficient in a product is a finite sum. Let $(\cdot)_+$ retain
nonnegative powers of $\partial$, and write $P_3:=4(Q^3)_+$,
$P_5:=16(Q^5)_+$. Explicitly,
\begin{equation}\label{eq:laxoperators}
 \begin{aligned}
 P_3&=4\partial^3+6u\partial+3u_{-1},\\
 P_5&=16\partial^5+40u\partial^3+60u_{-1}\partial^2
            +(50u_{-2}+30u^2)\partial+15u_{-3}+30uu_{-1}.
 \end{aligned}
\end{equation}
Their commutators with $L$ have order zero, with coefficients exactly
$u_{-3}+6uu_{-1}$ and $u_{-5}+10uu_{-3}+20u_{-1}u_{-2}+30u^2u_{-1}$. Hence
\begin{equation}\label{eq:lax}
 X_0L=[P_3,L],\qquad Y_0L=[P_5,L].
\end{equation}
This identifies the normalizations of the commuting classical flows.
The Lax identities extend to $Q^{2k-1}$. Taking residues gives
conservation because the residue of a pseudodifferential commutator
is a total derivative. This standard trace identity follows directly
from \eqref{eq:pdomult}: for $r=m+n+1\ge0$ the residue of
$[a\partial^m,b\partial^n]$ is
$\binom mr a\partial^r b-\binom nr b\partial^r a$, which is zero
modulo total derivatives since $\binom nr=(-1)^r\binom mr$.
Finally, setting all positive-order jets to zero gives
\[
 \ev\Res_\partial L^{k-1/2}
   =\binom{k-1/2}{k}u^k=\frac{\binom{2k}{k}}{4^k}u^k.
\]
Assigning weight one to $\partial$ makes $Q$ homogeneous of weight
one; its $\partial^{-1}$ coefficient in $Q^{2k-1}$ therefore has
weight $2k$. Reduction preserves this weight and the derivative-free
term, proving the assertions.
\end{proof}

\subsection{Kernel after reduction and extraction}
Let $D_3:=\sum_{j\ge0}u_{-j-3}\partial_{u_{-j}}$, so $X_0=\XH+D_3$.
Composing $X_0$ first with reduction $\Pi$ and then with $\extract$
gives
\begin{equation}\label{eq:C0}
 M_0:=\extract\Pi(N+D_3),\qquad
 \extract\Pi X_0=1+M_0-\ev.
\end{equation}
The order in Proposition~\ref{prop:finite} makes $M_0$ nilpotent:
$D_3$ preserves factor number, so $\extract\Pi D_3$ lowers it.
Equation~\eqref{eq:hlimit} also identifies $M_0$ with $M|_{\hh=0}$.

\begin{lemma}[Kernel of the extracted classical action]\label{lem:kernelBX}
For every positive weight $w$,
\begin{equation}\label{eq:kernelBX}
 \ker(\extract\Pi X_0|_{\Fcal_{0,w}})=
 \begin{cases}\Q I_k,&w=2k,\\0,&w\text{ odd}.\end{cases}
\end{equation}
The even-weight line is contained in $\ker X_0\cap\ker Y_0$.
In particular, $X_0$ is injective on positive odd-weight spaces.
\end{lemma}
\begin{proof}
An element of this kernel satisfies
$g=(1+M_0)^{-1}\ev g$. If $w$ is odd, $\ev g=0$ because every
power of $u$ has even weight; hence $g=0$. If $w=2k$, $\ev g$ is
a scalar multiple of $u^k$, so the kernel has dimension at most one.
The normalized class $I_k$ of Lemma~\ref{lem:rho} belongs to this
kernel and is killed by both flows, proving equality.
\end{proof}

\subsection{Middle exactness of the KdV Koszul complex}
The commuting endomorphisms $X_0$ and $Y_0$ define a Koszul complex
on the graded vector space $\Fcal_0$. Its middle exactness says that
every compatible pair $(a,b)$ comes from a common density. The proof
reduces this assertion, one factor number at a time, to divisibility
in a polynomial ring.

\begin{lemma}[Middle exactness]\label{lem:pair}
For even $d=2k\ge4$, the sequence
\begin{equation}\label{eq:compatible}
 \Fcal_{0,d}\xrightarrow{(X_0,Y_0)}
 \Fcal_{0,d+3}\oplus\Fcal_{0,d+5}
 \xrightarrow{(Y_0,-X_0)}\Fcal_{0,d+8}
\end{equation}
is exact at the middle term. Equivalently, $Y_0a=X_0b$ implies
$(a,b)=(X_0g,Y_0g)$ for some $g\in\Fcal_{0,d}$.
\end{lemma}
\begin{proof}
\textbf{Step 1. Replace derivatives by symmetric polynomials.}
Translation preserves factor number, so $\Fcal_0$ is the direct sum
of its factor-homogeneous parts. For a product of $m$ jets, assign one
auxiliary variable $\xi_i$ to each factor and record its derivative
order as an exponent. Averaging over the assignments makes the result
independent of the order of the factors. This is the symbolic
representation of differential polynomials
\cite[Section~2.3, Definition~6]{MNW}:
\begin{equation}\label{eq:polarization}
 u_{-j_1}\cdots u_{-j_m}\longmapsto
 \frac1{m!}\sum_{\sigma\in S_m}
       \xi_{\sigma(1)}^{j_1}\cdots \xi_{\sigma(m)}^{j_m}.
\end{equation}
Applying $\partial$ differentiates each factor in turn, so it becomes
multiplication by $p_1=\xi_1+\cdots+\xi_m$. Passing to densities modulo
total derivatives therefore amounts to imposing $p_1=0$.
Consequently the factor-$m$ part of $\Fcal_0$ is
identified with the graded polynomial ring
\begin{equation}\label{eq:symbolring}
 \mathscr R_m:=\Q[\xi_1,\ldots,\xi_m]^{S_m}/(p_1)
             =\Q[e_2,\ldots,e_m],
\end{equation}
where $e_j$ is the elementary symmetric polynomial, of degree $j$.
Weight $w$ corresponds to symbol degree $w-2m$. The constant in
$\mathscr R_m$ represents $u^m$, not the vacuum. These rings are used
one factor number at a time, not as a multiplication on $\Fcal_0$.

The factor-preserving parts of $X_0,Y_0$ are $D_3,D_5$, where
$D_su_{-j}=u_{-j-s}$. Each $D_s$ raises the derivative order of one factor by $s$ and sums
over all factors. Under \eqref{eq:polarization} it therefore becomes
multiplication by the power sum $p_s=\sum_{i=1}^m\xi_i^s$:
\begin{equation}\label{eq:power}
 p_3=3e_3,\qquad p_5=5(e_5-e_2e_3).
\end{equation}
These are Newton's identities with $e_1=0$; set $e_j=0$ when $j>m$.
All nonlinear terms in $X_0,Y_0$ raise factor number.

By the observation after \eqref{eq:redexamples}, both odd-weight
residuals start at factor number at least three. In these coordinates,
this is also immediate from $\mathscr R_1=\Q$ and
$\mathscr R_2=\Q[e_2]$.

\textbf{Step 2. Divisibility.} Suppose the residual pair has no components
of factor number less than $m$. Let $a_m,b_m$ denote its factor-$m$
symbols, allowing either to be zero. Nonlinear terms cannot contribute
at this lowest factor number, so compatibility reads
\begin{equation}\label{eq:pairlowest}
 p_5a_m=p_3b_m.
\end{equation}
For $m=3,4$, the symbol $a_m$ has odd degree, and $e_3$ is the only
odd-degree generator of $\mathscr R_m$. Every monomial in $a_m$ is
therefore divisible by $e_3$. For $m\ge5$, reduce
\eqref{eq:pairlowest} modulo $e_3$. The result is
$e_5(a_m\bmod e_3)=0$ in the polynomial domain
$\mathscr R_m/(e_3)$. Again $e_3$ divides $a_m$.

\textbf{Step 3. Elimination.} Write $a_m=p_3q_m$. Cancellation in the domain $\mathscr R_m$ gives
$b_m=p_5q_m$. Its degree is
\[
 \deg q_m=(d+3-2m)-3=d-2m.
\]
Thus $q_m$ represents an $m$-factor source $g_m$ of weight $d$.
Subtract the \emph{full} pair $(X_0g_m,Y_0g_m)$. Its leading pieces remove
$(a_m,b_m)$, its remaining pieces have higher factor number, and
$[X_0,Y_0]=0$ preserves compatibility. No earlier component is recreated.

\textbf{Step 4. Termination.} Perform this elimination for $m=3,\ldots,k$. Afterward the first
residual can have only factor number $k+1$, with symbol degree one;
there is no such element in \eqref{eq:symbolring}. Larger factor
numbers have negative symbol degree. Hence the first residual is zero.
Compatibility now says that $X_0$ kills the second residual, which
has odd weight. Lemma~\ref{lem:kernelBX} makes it zero as well.
The sum of the $g_m$ is the required $g$. For $k=2$ the elimination
range is empty and the same final argument applies.
\end{proof}

For $m\ge5$, divisibility uses that $e_5$ is a non-zero-divisor
modulo $e_3$; for $m=3,4$, it uses parity instead. This is why the
argument proves middle exactness in the stated weights without
requiring a regular sequence at every factor number.

\subsection{Compatibility and extraction determine a pair}
We now combine middle exactness with the common conserved line.
Compatibility represents a pair as $(X_0g,Y_0g)$; the additional
condition $\extract\Pi X_0g=0$ then puts $g$ on that line and makes
both entries vanish.
\begin{proposition}\label{prop:test}
For every even $d\ge4$, the following map is injective:
\begin{equation}\label{eq:test0}
 \begin{aligned}
 \Theta_0:\Fcal_{0,d+3}\oplus\Fcal_{0,d+5}
        &\longrightarrow\Fcal_{0,d+8}\oplus\Fcal_{0,d},\\
 (a,b)&\longmapsto(Y_0a-X_0b,\,\extract a)
 \end{aligned}
\end{equation}
\end{proposition}
\begin{proof}
For a pair in its kernel, Lemma~\ref{lem:pair} gives
$(a,b)=(X_0g,Y_0g)$. The second equation is $\extract\Pi X_0g=0$.
Lemma~\ref{lem:kernelBX} makes $g$ a multiple of $I_k$, and both
$X_0$ and $Y_0$ annihilate that class. Thus $a=b=0$.
\end{proof}

\section{Conservation, uniqueness, and commutativity}\label{sec:quantum}
Proposition~\ref{prop:projectedconservation} gives
$\extract\Pi X_\hh F_k=0$. To obtain conservation, we must deduce
$\Pi X_\hh F_k=0$. The classical limits are now identified by
\eqref{eq:hlimit}; we lift the injectivity of
Proposition~\ref{prop:test} to formal power series in $\hh$ and apply
it to the compatible pair $(X_\hh F_k,Y_\hh F_k)$. The finite inverse
then proves uniqueness, and the Jacobi identity gives commutativity.
\subsection{Lifting the injectivity to formal power series}
Extend scalars at each fixed weight:
\[
 \widehat{\Fcal}_w:=\Fcal_w\otimes_{\base}\Q[[\hh]]
                  \simeq\Fcal_{0,w}[[\hh]].
\]
Write $\widehat{\Fcal}=\bigoplus_w\widehat{\Fcal}_w$. These are
finite-rank free modules at each fixed weight. For even $d\ge4$, on
$\widehat{\Fcal}_{d+3}\oplus\widehat{\Fcal}_{d+5}$ define
\begin{equation}\label{eq:testh}
 \Theta_\hh(a,b):=\bigl(Y_\hh a-X_\hh b,\,\extract a\bigr).
\end{equation}
Its constant term is the map $\Theta_0$ of Proposition~\ref{prop:test}.

\begin{lemma}[Injectivity over formal power series]\label{lem:formaltest}
The map $\Theta_\hh$ is injective. For every $k\ge2$ and every
$f\in\widehat{\Fcal}_{2k}$,
\begin{equation}\label{eq:implication}
 \extract\Pi X_\hh f=0\quad\Longrightarrow\quad
                  X_\hh f=Y_\hh f=0.
\end{equation}
\end{lemma}
\begin{proof}
If a nonzero formal pair lies in $\ker\Theta_\hh$, its first nonzero
coefficient in $\hh$ lies in $\ker\Theta_0$, contrary to
Proposition~\ref{prop:test}. Hence $\Theta_\hh$ is injective. Since
$[X_\hh,Y_\hh]=0$ by \eqref{eq:Ydef}, we have
\[
 \Theta_\hh(X_\hh f,Y_\hh f)=(0,\extract\Pi X_\hh f).
\]
The stated implication follows.
\end{proof}

\subsection{Conservation and the classical limit of the constructed densities}
For $k\ge2$, the polynomial $F_k$ has coefficients in $\Q[\hh]$ and
satisfies $\extract\Pi X_\hh F_k=0$ by \eqref{eq:projectedeq}.
Applying Lemma~\ref{lem:formaltest} over $\Q[[\hh]]$ gives
\begin{equation}\label{eq:conservationproved}
 \boxed{X_\hh F_k=Y_\hh F_k=0\quad\text{in }\widehat{\Fcal}.}
\end{equation}
Because $\Fcal$ is a free $\Q[\hh]$-module with a fixed reduced
basis, the equality already holds over $\Q[\hh]$. In VOA notation
it is exactly
$\hh^{-1}H_4(0)H_{2k}\in T\Vcal$,
which is the required conservation statement.

Taking the constant term gives $X_0F_k(0)=0$, together with
$\ev F_k(0)=u^k$. Lemma~\ref{lem:kernelBX}, independently derived
from the Lax residues, forces
\begin{equation}\label{eq:classicallimitproved}
 \boxed{F_k(0)=I_k.}
\end{equation}
This proves the classical-limit part of the theorem.

\subsection{Uniqueness and commutativity from the Jacobi identity}
Let $g\in\Fcal_w$ be a homogeneous class annihilated by $X_\hh$.
The identity \eqref{eq:projectedoperator} and nilpotence of $M$ give
\begin{equation}\label{eq:uniqueness}
 (1+M)g=\ev g,\qquad
 \boxed{g=(1+M)^{-1}\ev g.}
\end{equation}
Derivative-free monomials have even positive weight, so this equation
forces $g=0$ when $w$ is odd. At weight $2k$ it forces $g$ to be a
multiple of $F_k$; conservation gives the required nonzero class.
This includes $F_1=u$, whose centrality was checked above.
These kernel conclusions hold over $\Q[\hh]$ and at every
specialization $\hh\in\C$. At $\hh=0$ they concern $X_0$,
rather than the zero-mode bracket, which vanishes in that specialization.

Write $h_k:=[H_{2k}]$ and $h_2:=[H_4]$ in the Lie algebra
$\Vcal/T\Vcal$, with bracket $\{[a],[b]\}=[a(0)b]$.
Conservation gives $\{h_2,h_k\}=0$. By the Jacobi identity,
\[
 \{h_2,\{h_k,h_l\}\}=\{\{h_2,h_k\},h_l\}
             +\{h_k,\{h_2,h_l\}\}=0.
\]
In the reduced model, $\{h_2,a\}=\hh X_\hh a$. Because the local-functional
module is $\hh$-torsion free, the Jacobi identity implies
$X_\hh\{h_k,h_l\}=0$. This bracket has weight $2k+2l-1$,
which is odd; it therefore vanishes by the preceding kernel theorem. Thus
\begin{equation}\label{eq:commutationproved}
 \boxed{H_{2k}(0)H_{2l}\in T\Vcal.}
\end{equation}
On every module the corresponding zero-mode operators commute,
by $[a(0),b(0)]=(a(0)b)(0)$ and $(Ta)(0)=0$.

\subsection{Specialization to every central charge}
The preceding argument uses the rescaling $c=6/\hh$. To include
$c=0$, we return to the original Virasoro coordinates and show that
both the constructed densities and the inverse identity are polynomial
in $c$. This also proves the uniqueness statement at $c=0$.
To distinguish the two PBW coordinate systems, write $v_{-j}$ for the
coordinate of $T^j\omega$ and $v:=v_{0}$. The map $\sigma_\omega$
is defined by \eqref{eq:sym} with $u_{-j}$ replaced by $v_{-j}$ and
$U_{-j}$ replaced by $j!L_{-j-2}$. Reduction, evaluation, and extraction
use the same formulas in the $v$-variables; we retain the symbols
$\Pi,\ev,\extract$ on this polynomial space.

After adjoining $\hh^{-1}$, define $S_\hh$ by
$S_\hh f(u_{0},u_{-1},\ldots):=f(\hh v_{0},\hh v_{-1},\ldots)$.
Thus it multiplies an $m$-factor monomial by $\hh^m$ and changes
its coordinate letters. Set
\begin{equation}\label{eq:reverse}
 \boxed{F_k^{(c)}:=\left.\hh^{-k}S_\hh F_k\right|_{\hh=6/c}.}
\end{equation}
For an $m$-factor monomial $\mon$, Proposition~\ref{prop:finite}
allows only terms $\hh^q\mon$ with $0\le q\le k-m$. Under
\eqref{eq:reverse}, such a term becomes
\[
 (c/6)^{k-m-q}\mon(v_{0},v_{-1},\ldots).
\]
Consequently $F_k^{(c)}\in\Q[c][v_{0},v_{-1},\ldots]$, its $m$-factor
coefficients have $c$-degree at most $k-m$, and
$\ev F_k^{(c)}=v^k$. Equation~\eqref{eq:reverse} therefore defines
the density at $c=0$ by polynomial evaluation.

Let the original weight-four action be
\[
 A_c:=\sigma_\omega^{-1}(\omega(-1)\omega)(0)\sigma_\omega.
\]
The PBW identifications satisfy $\sigma=\sigma_\omega S_\hh$ after
$c=6/\hh$. Since $H_4=\hh^2\omega(-1)\omega$, we have
\begin{equation}\label{eq:operatorscale}
 X_\hh=\hh S_\hh^{-1}A_{6/\hh}S_\hh.
\end{equation}
Likewise $\sigma_\omega(F_k^{(c)})=\hh^{-k}H_{2k}$ after localization.
The conservation and commutativity identities thus hold after the
substitution $c=6/\hh$. Virasoro PBW reordering makes the original
mode action polynomial in $c$, and the reduced basis is free over
$\Q[c]$. Injectivity of $\Q[c]\hookrightarrow\Q((\hh))$ then proves
these identities over $\Q[c]$, including their evaluation at $c=0$.

For the centralizer statement, put
$M_c:=\left.S_\hh M S_\hh^{-1}\right|_{\hh=6/c}$.
A term of $M$ lowering factor number by $r$ has parameter degree
at most $r$. Conjugation multiplies it by $\hh^{-r}$, so $M_c$
also has polynomial coefficients in $c$. Both $\Pi$ and $\ev$
commute with $S_\hh$, while
$\extract S_\hh=\hh S_\hh\extract$. Conjugating
\eqref{eq:projectedoperator} and using \eqref{eq:operatorscale} gives
\begin{equation}\label{eq:centralchargeidentity}
 \boxed{\extract\Pi A_c=1+M_c-\ev.}
\end{equation}
The decreasing order of Proposition~\ref{prop:finite} is unchanged,
so $M_c$ is nilpotent after every specialization. If $A_cg=0$ in
the local-functional space, then
$g=(1+M_c)^{-1}\ev g$. Odd weights therefore give zero, and an
even-weight kernel has dimension at most one. The conserved class
$F_k^{(c)}$, whose evaluation is $v^k$, is nonzero at every $c$
and spans that line. This completes the proof of
Theorem~\ref{thm:main}.

\clearpage
\appendix
\section{Low-weight formulas and exact checks}\label{app:computations}
We first display the vectors $H^{(c)}_{2k}$ in the original Virasoro
generators $L_n$ and conformal vector $\omega$, at central charge $c$.
The later calculations retain the rescaled coordinates used in the
construction. These exact examples illustrate the formulas; the proofs
in the main text establish the results at all weights.

\subsection{Conserved vectors in the original Virasoro notation}
Recall $H_{2k}^{(c)}=\sigma_\omega(F_k^{(c)})$ from
Theorem~\ref{thm:main}, with $H_2^{(c)}=\omega=L_{-2}\vac$.
For the first nontrivial weights, the literal symmetric representatives are
\begin{equation}\label{eq:Hlow}
\begin{aligned}
 H_4^{(c)}&=L_{-2}^2\vac,\\
 H_6^{(c)}&=\left(L_{-2}^3-\frac{c+2}{12}L_{-3}^2\right)\vac,\\
 H_8^{(c)}&=\left[L_{-2}^4-\frac{c+8}{3}
       \left(L_{-3}^2L_{-2}+L_{-5}L_{-3}+\frac43L_{-8}\right)
       +\frac{c^2+9c-61}{45}L_{-4}^2\right]\vac.
\end{aligned}
\end{equation}
Indeed, $v_{-j}$ corresponds to $T^j\omega=j!L_{-j-2}\vac$,
and the mixed symmetrization in the weight-eight formula is
\[
 \sigma_\omega(vv_{-1}^2)
 =\left(L_{-3}^2L_{-2}+L_{-5}L_{-3}+\frac43L_{-8}\right)\vac.
\]
The corresponding reduced differential polynomials are
\begin{equation}\label{eq:originalexamples}
 \begin{aligned}
 F_2^{(c)}&=v^2,\\
 F_3^{(c)}&=v^3-\frac{c+2}{12}v_{-1}^2,\\
 F_4^{(c)}&=v^4-\frac{c+8}{3}vv_{-1}^2
                   +\frac{c^2+9c-61}{180}v_{-2}^2.
 \end{aligned}
\end{equation}
Different ordered representatives can have different lower-factor
coefficients. For example, modulo $TV^c$,
\[
 \begin{aligned}
 H_6^{(c)}&\equiv
   \left(L_{-2}^3+\frac{c+2}{6}L_{-4}L_{-2}\right)\vac,\\
 H_8^{(c)}&\equiv
   \left(L_{-2}^4+\frac{c+8}{3}L_{-4}L_{-2}^2
      +\frac{2(c^2+14c-21)}{15}L_{-6}L_{-2}\right)\vac.
 \end{aligned}
\]
These congruences follow by applying $T=L_{-1}$ to shorter ordered
words. In particular, the difference between the two weight-six
representatives is
$-(c+2)T(L_{-3}L_{-2}\vac)/12$.
Comparison with cylinder charges also requires the coordinate change
discussed in the introduction; equality of the displayed vectors is
not asserted in a different coordinate convention.

\subsection{Explicit total-derivative identities}
The conservation identities can be verified before passing to the
quotient by displaying vectors $W_{k,l}^{(c)}$ such that
\[
 H_{2k}^{(c)}(0)H_{2l}^{(c)}=T W_{k,l}^{(c)}.
\]
Since $\omega(0)=T$, one may take $W_{1,l}^{(c)}=H_{2l}^{(c)}$.
For the weight-four pair, an explicit choice is
\begin{equation}\label{eq:W22}
 W_{2,2}^{(c)}=\left[
 4L_{-2}^3+\frac{2(c+5)}3L_{-4}L_{-2}
       +\frac{4-c}{6}L_{-3}^2+\frac{19c-22}{15}L_{-6}\right]\vac.
\end{equation}
This follows from \eqref{eq:Xu2} after inverse scaling and
symmetrization; it can also be checked directly using the Virasoro
commutator and $[T,L_{-n}]=(n-1)L_{-n-1}$.
For the distinct pair of weights four and six, one may take
\begin{equation}\label{eq:W23}
\begin{aligned}
 W_{2,3}^{(c)}=\Bigl[&\frac92L_{-2}^4
 +(c+11)L_{-4}L_{-2}^2-(c+2)L_{-3}^2L_{-2}
 +3(c+4)L_{-6}L_{-2}\\
 &-\frac{c^2+13c-14}{6}L_{-5}L_{-3}
 +\frac{c^2+c-29}{18}L_{-4}^2
 -\frac{19c^2-227c+226}{42}L_{-8}\Bigr]\vac.
\end{aligned}
\end{equation}
Substitution in the zero-mode sum verifies
$H_4^{(c)}(0)H_6^{(c)}=TW_{2,3}^{(c)}$.
Both identities use the literal representatives in \eqref{eq:Hlow}.

\subsection{Evaluating the finite inverse}
To exhibit the construction itself, return to the rescaled coordinates
$u_{-j}$ over $\Q[\hh]$. The reduced weight-eight polynomial is
\begin{equation}\label{eq:F4}
 \boxed{F_4=u^4-\left(2+\frac{8\hh}{3}\right)uu_{-1}^2
            +\left(\frac15+\frac{3\hh}{10}
                        -\frac{61\hh^2}{180}\right)u_{-2}^2.}
\end{equation}
Here the inverse in \eqref{eq:Fk} has just three nonzero summands:
\begin{equation}\label{eq:Mexample}
\begin{aligned}
 -Mu^4&=-\left(2+\frac{8\hh}{3}\right)uu_{-1}^2
           +\left(\frac{3\hh}{20}-\frac{11\hh^2}{60}\right)u_{-2}^2,\\
 M^2u^4&=\left(\frac15+\frac{3\hh}{20}
                           -\frac{7\hh^2}{45}\right)u_{-2}^2,
 \qquad M^3u^4=0.
\end{aligned}
\end{equation}
The calculation applies \eqref{eq:Xmodes}, reduction, and extraction
at each step. Inverse scaling by \eqref{eq:reverse} gives
$F_4^{(c)}$ in \eqref{eq:originalexamples}.

\subsection{The mode calculation at weight four}\label{app:modecheck}
We record the details of the calculation leading to \eqref{eq:Xu2}.
For the action on $u^2$, we use ordered products of the $U_{-j}$.
The terms with at least two factors are as follows; the omitted terms
in this table have at most one factor:
\[
 \begin{array}{c|l}
 r&\displaystyle \dfrac{2}{\hh(r-2)!}U_{-r+2}U_{r-1}U_{0}^2\vac\ \text{modulo one-factor vectors}\\\hline
 2&(4U_{-1}U_{0}^2+\hh U_{-3}U_{0})\vac\\
 3&8U_{-1}U_{0}^2\vac\\
 4&3\hh U_{-2}U_{-1}\vac\\
 5&(2+\tfrac83\hh)U_{-3}U_{0}\vac.
 \end{array}
\]
The additional nonzero $r=7$ contribution has only one factor.
The exact PBW conversions, including their lower-factor terms, are
\[
 \begin{aligned}
 \sigma(u^2u_{-1})&=(U_{-1}U_{0}^2+\tfrac\hh6 U_{-3}U_{0}+\tfrac{\hh^2}{120}U_{-5})\vac,\\
 \sigma(uu_{-3})&=(U_{-3}U_{0}+\tfrac{3\hh}{40}U_{-5})\vac,\\
 \sigma(u_{-1}u_{-2})&=(U_{-2}U_{-1}+\tfrac\hh{120}U_{-5})\vac,\qquad
 \sigma(u_{-5})=U_{-5}\vac.
 \end{aligned}
\]
They follow by reordering the two- and three-factor permutation averages.
Including the one-factor contributions, the full mode sum is
\[
 \left(12U_{-1}U_{0}^2+(2+\tfrac{11}3\hh)U_{-3}U_{0}
       +3\hh U_{-2}U_{-1}+(\tfrac3{10}\hh+\tfrac1{20}\hh^2)U_{-5}\right)\vac.
\]

Inverting these symmetric PBW identities gives \eqref{eq:Xu2}.

\section*{AI-assisted research and writing}
The problem studied in this paper was proposed by the author,
who set the research objectives and guided the scope, conventions,
and successive revisions. The work developed through sustained
interaction with OpenAI's ChatGPT (GPT-6 Pro, web interface) and
Codex (GPT-6 Astra). The author contributed mathematical
questions, choices of formulation, and critical revisions; the AI
systems suggested intermediate constructions, possible lines of
argument, and alternative ways to organize proofs. Some ideas were
developed from the author's proposals and corrections, while others
began as AI suggestions and were refined through subsequent exchanges.
AI assistance thus contributed to the development of mathematical
ideas and arguments as well as to their presentation.

ChatGPT and Codex were also used to draft and revise exposition,
write code for exact symbolic calculations, assist with checking
intermediate arguments and formulas, and locate literature for
verification. AI suggestions are not treated as independent validation
of mathematical claims. The conclusions rest on the proofs and cited
results, supplemented by exact computational checks. The author
retains responsibility for the final content.

\bigskip

\noindent Hongbo Zhao\par
\medskip
\noindent School of Mathematics and Statistics, Chongqing University, Chongqing, P.R.C. China\par
\noindent E-mail: \texttt{hongbozhao@cqu.edu.cn}

\end{document}